\documentclass[11pt]{article}
\usepackage{fancyhdr}
\usepackage{mathrsfs}
\usepackage{amsmath}
\usepackage{latexsym}
\usepackage{amssymb}
\usepackage{amsthm}
\usepackage{indentfirst}
\usepackage{titlesec}
\usepackage{titletoc}
\usepackage{graphicx}
\renewcommand{\baselinestretch}{1.1}

\begin{document}
%\begin{CJK*}{GBK}{song}
\title{ Classification of Element Systems
over Finite Commutative  Groups}
\author{ Junqin Li $^a$,
 Min Wu $^{a, b }$, \ \  Hengtai Wang $^{a}$, \ \ Shouchuan Zhang $^a$ \\
a: Department  of Mathematics, Hunan University, Changsha,
410082\\
b: Department of Mathematics, Tsinghua University, Beijing, 100084
 }
 \date{}

%=============================================== 自己的设置 ========================================================%

\newtheorem{Proposition}{\quad Proposition}[section]
\newtheorem{Theorem}[Proposition]{\quad Theorem}
\newtheorem{Definition}[Proposition]{\quad Definition}
\newtheorem{Corollary}[Proposition]{\quad Corollary}
\newtheorem{Lemma}[Proposition]{\quad Lemma}
\newtheorem{Example}[Proposition]{\quad Example}

\maketitle \addtocounter{section}{-1}

\numberwithin{equation}{section}

\begin {abstract} Directed graphs are widely applied in network.  
We  obtain the formula  computing  the number of isomorphic
classes of element systems with characters over finite commutative
group $G$.


\end {abstract}

\numberwithin{equation}{section}
%%\numberwithin{eqnarray}{section}

%\pagenumbering{arabic}

\section{Introduction}\label {s0}

Classification of  Hopf algebras was developed and popularized in
the last decade of the twentieth century,  which would have
applications to a number of other areas of mathematics, aside from
its intrinsic algebraic interest. In mathematical physics,
Drinfeld's and Jambo's work was to provide solutions to quantum
Yang-Baxter equation. In conformal field theory, I. Frenkel and Y.
Zhu have shown how to assign a Hopf algebra to any conformal field
theory model\cite {FZ92}. In topology, quasi-triangular and ribbon
Hopf algebras provide many invariants of knots, links, tangles and
3-manifolds\cite {He91, Ka97, Ra94, RT90}. In operator algebras,
Hopf algebras can be assigned as an invariant for certain
extensions.

Researches on the classification of Hopf algebra is in the
ascendant.  N. Andruskiewitsch and H. J. Schneider have obtained
interesting result in   classification of finite-dimensional
pointed Hopf algebras
 with commutative coradical  \cite {AS98a, AS98b, AS02, AS00}.
 More recently, they  have also researched this problem in case of
 non-commutative coradical. Pavel Etingof and Shlomo Gelaki gave the complete
 and explicit classification of finite-dimensional
triangular Hopf algebras over an algebraically closed field $k$ of
characteristic $0$ \cite {EG00}. The classification of monomial
Hopf algebras, which are a class of co-path Hopf algebras, and
simple-pointed sub-Hopf algebras of co-path Hopf algebras were
recently obtained in \cite {CHYZ04} and \cite {OZ04},
respectively.

Assume that $k$ is an algebraically closed field of characteristic
zero with a primitive $\mid \! G \! \mid $th root of 1 and $G$ is
a finite abelian group.
 Element systems with characters can be
applied to classify quiver Hopf algebras, multiple Taft algebras
over $G$ and Nichols algebras in $^{FG}_{FG} {\cal YD}$ (see
\cite[Theorem 3, Theorem 4]{ZZC04} ).
 In this paper, we  obtained the formula  computing  the number of isomorphic classes
of element systems with characters over finite commutative group
$G$.

In  Section \ref {s1} we give the explicit formula computing the
number of isomorphic classes of element systems over finite cycle
$p$-groups. In  Section \ref{s2} we give the explicit formula
computing the number of isomorphic classes of element systems over
finite cycle groups. In Section \ref {s3} we give the explicit
formula computing the number of isomorphic classes of element
systems of primary commutative groups. In Section \ref {s4} we
give the formula computing the number of isomorphic classes of
element systems over finite commutative groups. Unfortunately, we
have not found the explicit formula in the general case; we only
change the problem into the numbers of  solutions of congruence
class equations  (\ref {e4.3}) and (\ref {e4.4}).

\section {Preliminaries} \label {s1'}

Unless specified otherwise, in the paper we have the following
assumption and notations. $G$ is a abelian group with order $m$;
$F$ is a field containing a primitive $m$th root of 1;Aut$G$ and
Inn$G$ denote the automorphism and inner automorphism group
respectively; $1$ denotes the unity element of $G;$ $\widehat{G}$
denotes the set of characters $G$, where a character of $G$ is   a
group homomorphism from $G$ to $F -\{0\}$.
 $C_m$ denotes a cycle group with order $m$; $\mathbb{ Z}$
denotes the  set of all  integers; $\mathbb{N}$ denotes the set of
all positive integers;  $\mathbb{Z}_m$ denotes the ring of integers
modulo $m$; $A^B$ denotes the cartesian product $\prod \limits_
{i\in B}A_i$, where $A_i = A$ for any $i \in B$. Since field $F$
contains a primitive $|G|$th root of 1, $G$ has $|G|$ elements, i.e.
$|\widehat{G}|=|G|. $

$\varphi(n)$ denotes the Euler function, i.e. $\varphi(n)$ is  the
number of elements  in set $\{ x \mid      (x,n )=1, 1 \le x \le n
\}$. For convenience, we denote $\bar i \in \mathbb{Z}_m$ by  $i$.
Obviously, $ \{i \mid (i, m)=1, 1\le i \le m\}$ is the set of all
invertible elements in multiplicative group $\mathbb{Z}_{m}$,
written $\mathbb{Z}_{m}^*$. It is clear that the number of elements
in set $\mathbb{Z}_{m}^* $ is $\varphi (m)$.

%对于 $\alpha$
%生成的$m$阶 循环群 $C_m = (\alpha )$,  Definition  $C_m ^* : = \{ \alpha ^i
%\mid  (i, m) = 1\}$. 易知 $Z_m^*$ 是 $C_m$的 $\varphi (m)$阶子群.

\begin{Lemma}\label{0.2} (i)  There are
$\frac{n!}{\lambda_1!\lambda_2!\cdots\lambda_n!1^{\lambda_1}2^{\lambda_2}\cdots
n^{\lambda_n}}$ permutations of type
$1^{\lambda_1}2^{\lambda_2}\cdots n^{\lambda_n}$ in $S_n$. (ii) If
$n=p^{e_1}_1p_2^{e_2}\cdots p_s^{e_s}$ and $p_1, p_2, \cdots, p_s $
are mutually different  prime numbers with positive integer $e_i$
for $i =1, 2, \cdots, s,$ then
$$
\varphi(n)=n\left(1-\frac{1}{p_1}\right)\left(1-\frac{1}{p_2}\right)\cdots\left(1-\frac{1}{p_s}\right).$$
In particular, $\varphi(p^{e_i})=p^{e_i}-p^{e_i-1}$.
  (iii) If $n\in
\mathbb{N},$ then $\sum\limits_{d|n}\varphi(d)=n$. (iv) If
$\lambda_1, \lambda_2, \cdots, \lambda_n, n \in \mathbb{N}$, then
$\sum\limits_{\lambda_{1}+2\lambda_{2}+\cdots
+n\lambda_{n}=n}\frac{1}{\lambda_1!\lambda_2!\cdots\lambda_n!1^{\lambda_1}2^{\lambda_2}\cdots
n^{\lambda_n}} =1$.

\end{Lemma}
{\bf Proof. } (i) It follows from  \cite [Exercise  2.7.7] {Hu05}.

(ii) and  (iii) follow from \cite [ Theorem  3.3.1,  Theorem 3.3.1]
{BB99}.

(iv) Considering  $$ \sum\limits_{\lambda_{1}+2\lambda_{2}+\cdots
+n\lambda_{n}=n}\frac{n!}{\lambda_1!\lambda_2!\cdots\lambda_n!1^{\lambda_1}2^{\lambda_2}\cdots
n^{\lambda_n}}=n!,  $$ we complete the proof. $\Box$

A group $G$ is said  to act on a non-empty set $\Omega$, if there is
a map  $G\times \Omega\rightarrow\Omega$, denoted by $(g, x)\mapsto
g\circ x$, such that for all $x\in \Omega$ and $g_1, g_2\in G$:
$$
(g_1g_2)\circ x=g_1\circ(g_2\circ x)\qquad and\qquad   1\circ x=x,$$
where 1 denote the unity element of $G.$

For each $x\in \Omega$, let
$$ G_x:=\{g\circ x \mid  g\in G\},$$
 called the orbit of $G$ on $\Omega$. For each $g\in G$, let
 $$ F_g:=\{x\in \Omega \mid  g\circ x=x\},$$
 called the fixed point set of  $g$.

Burnside's lemma, sometimes also called Burnside's counting theorem,
which is  useful to compute the number of orbits.

\begin{Theorem}\label{0.6} (See \cite [Theorem 2.9.3.1]{Hu05})
( Burnside's Lemma )  Let $G$ be a finite group that acts on  a
finite set $\Omega$, then the number $\mathcal {N}$ of orbits is
given by the following formula:
$$
\mathcal {N}=\frac{1}{ \mid  G \mid  }\sum\limits_{g\in G} \mid F_g
\mid  . $$

\end {Theorem}

\begin {Definition}\label{0.8}

$(G, \overrightarrow{g}, \overrightarrow{\chi}, J)$ is called an
 element system with characters  (simply, ESC) if $G$ is a
group, $J$ is a set, $\overrightarrow{ g } = \{g_i\} _{ i\in J} \in
Z(G)^J$ and $\overrightarrow{ \chi }= \{\chi_i\}_{ i \in J} \in
\widehat{G}^J $ with $ g_i \in Z(G)$ and $\chi _i \in \widehat G$.
$ESC (G, \overrightarrow{g}, \overrightarrow{\chi}, J)$ and $ESC
(G', \overrightarrow{g'}, \overrightarrow{\chi'}, J')$ are said to
be isomorphic if there exist a group isomorphism $\phi: G
\rightarrow G'$ and a bijective map $\sigma: J\rightarrow J'$ such
that $\phi(g_i)=g'_{\sigma(i)}$ and $\chi'_{\sigma(i)}\phi=\chi_i$
for any $i \in J$.
\end {Definition}

$ ESC(G, \overrightarrow{g}, \overrightarrow{\chi}, J)$    can be
written as  $ ESC(G, g_i, \chi _i; i\in J ) $ for convenience.

Given a finite commutative group and a positive integer $n$,
define \ \ \ \ \ \ $\Omega (G, n) : = $ $\{ (G,
\overrightarrow{g}, \overrightarrow{\chi}, J) \mid (G,
\overrightarrow{g}, \overrightarrow{\chi}, J) $ $ \hbox {  is an }
ESC
 \hbox { and  } $ $ J = \{1, 2, \cdots, n\}\}$. Let $\left(\begin{array}{c} j\\
g_{j}\\
\chi_{j} \end{array}\right)$ denote   $ESC(G,g_{j},\chi_{j}; j\in
 J) $ $\in \Omega (G, n)$ in short.
 Define the action $\circ$ of group $M:=$Aut$G\times S_n$  on $\Omega (G,
n)$ as follows:
\begin {eqnarray}\label {e0.1}
 (\phi,\sigma)\circ \left(\begin{array}{c} j\\
g_{j}\\
\chi_{j} \end{array}\right)
=\left(\begin{array}{c} \sigma(j)\\
\phi(g_{j})\\
\chi_{j}\phi^{-1} \end{array}\right).
\end {eqnarray}
It is clear that it is an action. Indeed,
$$
(\phi,\sigma)\circ\left((\phi',\sigma')\circ \left(\begin{array}{c} j\\
g_{j}\\
\chi_{j} \end{array}\right) \right)=(\phi,\sigma)\circ\left(\begin{array}{c} \sigma'(j)\\
\phi'(g_{j})\\
\chi_{j}\phi^{\prime-1} \end{array}\right)=\left(\begin{array}{c} \sigma\sigma'(j)\\
\phi\phi'(g_{j})\\
\chi_{j}\phi^{\prime-1}\phi^{-1 }
\end{array}\right)=(\phi\phi',\sigma\sigma')\circ\left(\begin{array}{c} j\\
g_{j}\\
\chi_{j} \end{array}\right),
$$ and
$$
(1,1)\circ\left(\begin{array}{c} j\\
g_{j}\\
\chi_{j} \end{array}\right)=\left(\begin{array}{c} j\\
g_{j}\\
\chi_{j} \end{array}\right).
$$
Let ${\cal N}(G, n)$ denote the number of isomorphic classes in $
\Omega (G, r)$. This article is mostly devoted to investigate the
formula of ${\cal N}(G, n)$. It is clear that each orbit of
$\Omega(G, n)$ represents an isomorphic class of ESC's. As a
result, ${\cal N}(G, n)$ is equal to the number of orbits in
$\Omega(G, n )$.

\section{Cycle $p$-group}\label {s1}

In this section we give the formula computing the number of
isomorphic classes of element systems over finite cycle $p$-groups.

Let $G=\langle g\rangle$ and  $\mid \! G \! \mid = m = p^e.$ It is
clear that $\hat{G}$ consists of the following maps:
\begin{align*}\chi^{(l)}:G &\rightarrow F
\\g &\mapsto \omega^{l},
\end{align*} where $\omega$ is a primitive $m$th root of 1;
Aut$G$ consists of the following maps:
\begin{align*}
\phi_{i}:G &\rightarrow G\\
g &\mapsto g^{i},
\end{align*}where $(i,m)=1$.
%若$C^{\ast}_h$表示$Z_h$中的单位构成的子群.则Aut$G\cong
%C^*_{m}$且
That is, $\mathrm{Aut}G \cong \mathbb{Z}_m^*$ (See \cite [Theorem
2.3.3] {Zh82}). Thus $|\mathrm{Aut}G|=|Z^*_m|=\varphi(m)$.

 Since $\chi^{(l_j)}\phi_i^{-1}(g)=
\chi^{(l_j)}\phi_{i^{-1}}(g)=\chi^{(l_j)}(g^{i^{-1}})=\omega^{l_ji^{-1}}$,
$\chi^{(l_j)}\phi_i^{-1}=\chi^{(l_{j}i^{-1})}$. By (\ref{e0.1}),
we have
$$
 (\phi_{i},\sigma)\circ\left(\begin{array}{c} j\\
g^{k_j}\\
\chi^{(l_j)} \end{array}\right)
 =\left(\begin{array}{c} \sigma(j)\\
\phi_i(g^{k_j})\\
\chi^{(l_j)}\phi_i^{-1} \end{array}\right)  =\left(\begin{array}{c} \sigma(j)\\
g^{i\cdot
k_{j}}\\
\chi^{(l_{j}i^{-1})}\end{array}\right),$$ where $i^{-1}$ is the
inverse of $i$ in group  $\mathbb{Z}_m^*$.

\begin{Definition}\label{1.1}
 Assume $(i,p)=1$ and $p$ is a prime number. Define  $\vec{\delta}_e(i)$ $:= $ $(\delta_1
(i),\delta_2(i), $ $ \ldots,\delta_e(i))$, where $\delta_s(i)$ is
the order of $i$ in $Z^{*}_{p^s}$ for  $1\le s\le e$.
\end{Definition}

{\bf Remark.} It is easy to check
$\vec{\delta}_e(i)=\vec{\delta}_e(i^{-1})$.

\begin{Lemma}\label{1.3}
 If $p$ is an odd prime number, then there exists $\alpha\in
 \mathbb{Z}$ such that
$$Z^{*}_{p^s}=\langle\alpha\rangle=\{1,\alpha,\ldots,\alpha^{\varphi(p^s)-1}\}$$
 for any positive integer $s$. That is,
$\alpha$ is a common generator of $Z^{*}_{p^s}(s=1,2,\ldots)$.

If $p=2$, then there exist $\beta\in \mathbb{Z}$(for example, $\beta
=5$ ) such that
$$Z^{*}_{p^s}=\langle-1\rangle\times\langle\beta\rangle=\{1,-1\}\times\{1,\beta,\ldots,\beta^{2^{s-2}-1}\}$$
when $s>2$.

\end{Lemma}
{\bf Proof.} It follows from \cite[  Theorem 5.2.3;  Theorem
5.2.1]{BB99}. $\Box$

\begin{Lemma}\label{1.4}
(i) If  $d \mid m $, then  the number of elements with order $d$ in
cycle group $(\mathbb{Z}_m , +)$ is exactly  $\varphi (d)$.

(ii) If $p$ is an odd prime number or $e\leq 2$, then
$\vec{\delta}_e(i)$ is a case of  following for any $i\in
\mathbb{Z}_{m}^*$:
$$
\vec{\delta}_e^{(k,d)}=(d,\ldots,\overset{k}{d},pd,p^2d,\ldots,p^{e-k}d),\
d|p-1,k=1,\ldots,e.$$ Conversely, if $d|p-1$ and $1\leq k \leq e$,
then there exactly exist $\varphi(p^{e-k}d)$ elements in
$Z^{*}_{p^e}$ satisfying $\vec{\delta}_e(i)=\vec{\delta}^{(k,d)}_e$.

(iii) If $p=2$ with $e\geq3$, then $\vec{\delta}_e(i)$ is a case of
the following  for any $i\in \mathbb{Z}_{m}^*$:
$$
\vec{\delta}_e^{(k,d)}=(1,d,\ldots,\overset{k}{d},2,2^2,\ldots,2^{e-k}
),\ d=1 \mbox{ or } 2, k=2, \ldots, e.$$ Conversely, if $d=1$ or 2
and $2\leq k \leq e$, then there exactly exist $\varphi(2^{e-k})$
elements in  $Z^{*}_{2^e}$ satisfying $\vec{\delta}_e
(i)=\vec{\delta}_e^{(k,d)}$.
\end{Lemma}

{\bf Remark.} In (iii) above,
$\vec{\delta}_e^{(e-1,2)}=\vec{\delta}_e^{(e,2)}=(1,2,\ldots,2)$,
$\vec{\delta}_e^{(e,1)}=(1,1,\ldots,1)$.

{\bf Proof.} (i) Let  $A =\{i \mid 1\le i \le d, (i, d) =1\}$ and $B
=\{i \mid 1\le i \le m,$ the order of $i$ in  $( \mathbb{Z}_m, +) $
is $ d\}$.  Define a map $\phi: \ A \rightarrow B $ such that  $\phi
(i) = i \frac {m} {d}$ for any  $i \in A.$ It is clear that $\phi$
is bijective.

(ii) By  Lemma\ref{1.3},  we can assume $i=\alpha^{\nu}$ and
$\nu=p^{k-1}c$ with $e\ge k\geq 1$ and $(p,c)=1.$ For  any $1\leq
s\leq e$, since the order of $\alpha$ in $Z^{*}_{p^s}$ is
$\varphi(p^s)=p^{s-1}(p-1)$, we have that the order of $i$ in
$Z^{*}_{p^s}$ is
\begin{equation}\label{e1.2}
\delta_s
(i)=\frac{p^{s-1}(p-1)}{(\nu,p^{s-1}(p-1))}=\frac{p^{s-1}(p-1)}{(p^{k-1}c,p^{s-1}(p-1))}=\left\{
\begin{array}{ll}
p^{s-k}\frac{p-1}{(c,p-1)} & \hbox{ when $s\geq k$}\\
\frac{p-1}{(c,p-1)} & \hbox{when $s<k$}.
\end{array}
\right.
\end{equation}
Let $\frac{p-1}{(c,p-1)}=d$. Obviously $d|p-1$. To complete the
proof of (ii), we now  show the following conclusion:

If  $i, i' \in \mathbb{Z}_{p^e}$ and  $(ii',p)=1, $ then
$\vec{\delta}_e(i)=\vec{\delta}_e(i')$ if and only if $i$ and $i'$
have the same orders in $Z^{*}_{p^e}$. Indeed, the necessity is
obvious. Now we show the sufficiency. By  Lemma\ref{1.3}, we can
assume $i=\alpha^\nu$ and $ i'=\alpha^{\nu'}$ with $\nu=p^{k-1}c$,
$\nu'=p^{k'-1}c'$,  $(cc',p)=1$. Since $i$ and $i'$ have the same
orders in $Z^{*}_{p^e}$,  $k=k',(c,p-1)=(c', p-1)$ by (\ref{e1.2}).
Applying  (\ref {e1.2}), we have
$\vec{\delta}_e(i)=\vec{\delta}_e(i')$.

Next we come back to show (ii). For any $d$ and $k$ with $d|p-1$
and $1\leq k \leq e$, Let $i=\alpha^{p^{k-1}\frac{(p-1)}{d}}$. It
is easy to check $\vec{\delta}_e(i)=\vec{\delta}_e^{(k,d)}$.
Considering that there exactly exist $\varphi(p^{e-k}d)$ elements
whose orders are $p^{e-k}d$  in $Z^{*}_{p^e}$, we complete the
proof of (ii).

 (iii) If $i=\beta^\nu$, we can similarly show the first part of (iii). Notice $\delta_1 (i)=1.$
 If $i=-\beta^{\nu}$ and $\nu=2^{k-2}c$ with
 $(c,2)=1$ and  $e \ge k \ge 2$, then  the  order of  $i$ in
 $Z^{*}_{2^s}$ is
$$\delta_s(i)=[2, \frac{2^{s-2}}{(\nu, 2^{s-2})}]=[2, \frac{2^{s-2}}{(2^{k-2}c,
 2^{s-2})}]=\left\{
\begin{array}{ll}
2^{s-k} & \hbox{when } s> k \\
2 & \hbox{when $s\leq k$}
\end{array}
\right.
$$  when $s\ge 2$.
Obviously, $\delta_1(i)=1$ and
$\vec{\delta}_e(i)=\vec{\delta}_e^{(k,2)}$. It is easy to show the
second part of (iii) by meas of the method similar to proof of (ii).
$\Box$

 Next we compute the fixed point set $F_{(\phi_i,\sigma)}$ of $(\phi_i,\sigma)$.
Assume $\vec{\delta}_e(i)=\vec{\delta}_e^{(k,d)}$ and that
 the type of  $\sigma\in S_n$ is  $1^{\lambda_{1}}2^{\lambda_{2}}\cdots
n^{\lambda_{n}}$.
If $\left(\begin{array}{c} j\\
g^{k_{j}}\\
\omega ^{l_{j}} \end{array}\right)\in F_{(\phi_i,\sigma)},$ then
$$\left(\begin{array}{c} \sigma(j)\\
g^{ik_{j}}\\
\omega ^{i^{-1}l_{j}} \end{array}\right)=\left(\begin{array}{c} j\\
g^{k_{j}}\\
\omega ^{l_{j}} \end{array}\right)=\left(\begin{array}{c} \sigma(j)\\
g^{k_{\sigma(j)}}\\
\omega ^{l_{\sigma(j)}}
\end{array}\right),$$ which implies
\begin{equation}\label{e1.3}
ik_{j}\equiv k_{\sigma(j)}(\mathrm{mod}\ p^e),
\end{equation}
and
\begin{equation}\label{e1.4}
i^{-1}l_{j}\equiv l_{\sigma(j)}(\mathrm{mod}\ p^e)
\end{equation} for $1\leq j\leq n$.

It is clear that  (\ref {e1.3} ) and (\ref {e1.4} ) are independent
each other, and they hold if and only if every cycle, such as
$\tau=(j_{1}\,j_{2}\,\cdots j_{r})$, in independent cycle
decomposition of $\sigma$ satisfies the following two formulae:
\begin{equation*}\tag{\ref{e1.3}$'$}
ik_{j}\equiv k_{\tau(j)}(\mathrm{mod}\ p^e)
\end{equation*}
and
\begin{equation*}\tag{\ref{e1.4}$'$}
i^{-1}l_{j}\equiv l_{\tau(j)}(\mathrm{mod}\ p^e).
\end{equation*} for $j=j_1,\ldots,j_r$.

For (\ref{e1.3}$'$), we have  \setcounter{equation}{4}
\begin{eqnarray}\label{e1.5}
k_{j_{2}}&\equiv &ik_{j_{1}}\equiv i^{1}k_{j_{1}}, \nonumber\\
k_{j_{3}}&\equiv &ik_{j_{2}}\equiv i^{2}k_{j_{1}},\nonumber\\
&\cdots & \nonumber\\
k_{j_{r}}&\equiv &ik_{j_{r-1}}\equiv i^{r-1}k_{j_{1}},\nonumber\\
k_{j_{1}}&\equiv &ik_{j_{r}}\equiv i^{r}k_{j_{1}}.
\end{eqnarray}
This implies the numbers of solutions  of both (\ref{e1.3}$'$) and
(\ref{e1.5}) are the same when we view $k_{j_1},\ldots,k_{j_r}$ as
indeterminates. (\ref{e1.5}) is equivalent to
\begin{equation}\tag{\ref{e1.5}$'$}
(i^r-1)k_{j_1}\equiv0(\mathrm{mod}\ p^e).
\end{equation}

Now we give the main result.
\begin{Theorem}\label{1.6}
Assume that  $G$ is a cycle group with order $p^e$ and  positive
integer $e$.

(i) If $p$ is an odd prime number or $e\leq2$, then
\begin{eqnarray}\label{e1.6}
{\cal N}(G,
n)&=&\frac{1}{p^{e-1}(p-1)}\sum_{d|p-1}\varphi(d)\sum_{\lambda_1+2\lambda_{2}+\cdots+n\lambda_n=n}
\frac{p^{2f_p(\lambda_1,\ldots,\lambda_n,\vec{\delta}_e^{(e,d)})}}{\lambda_1!\cdots\lambda_n!1^{\lambda_1}\cdots
n^{\lambda_n}}\nonumber\\
&&+\sum_{k=1}^{e-1}p^{-k}\sum_{d|p-1}\varphi(d)\sum_{\lambda_1+2\lambda_{2}+\cdots+n\lambda_n=n}
\frac{p^{2f_p(\lambda_1,\ldots,\lambda_n,\vec{\delta}_e^{(k,d)})}}{\lambda_1!\cdots\lambda_n!1^{\lambda_1}\cdots
n^{\lambda_n}}
\end{eqnarray}
\begin{eqnarray}\label{e1.6'}
&=& 1+\frac{1}{p^{e-1}(p-1)}\sum\limits_{d|p-1\atop d\leq
n}\varphi(d)\left(\sum\limits_{\lambda_{1}+2\lambda_{2}+\cdots
+n\lambda_{n}=n}\frac{p^{2e\sum\limits_{t=1}^{[n/d]}\lambda_{td}}}{\lambda_1!\cdots\lambda_n!1^{\lambda_1}\cdots
n^{\lambda_n}}-1\right) \nonumber\\
&+&\sum\limits_{k=1}^{e-1}p^{-k}\sum\limits_{d|p-1\atop d\leq n
}\varphi(d)\left(\sum\limits_{\lambda_{1}+2\lambda_{2}+\cdots
+n\lambda_{n}=n}\frac{p^{2f_p(\lambda_1,\ldots,\lambda_n,\vec{\delta}_e^{(k,d)})}}{\lambda_1!\cdots\lambda_n!1^{\lambda_1}\cdots
n^{\lambda_n}}-1\right),
\end{eqnarray}
where \begin{equation}\label{e1.8}
f_p(\lambda_1,\ldots,\lambda_n,\vec{\delta}_e^{(k,d)})=\sum\limits_{s=0}^{e-k-1}(k+s)\sum\limits_{1\leq
t\leq [n/(p^sd)]\atop (t,p)=1}
\lambda_{tp^sd}+e\sum\limits_{t=1}^{[n/p^{e-k}d]} \
\lambda_{tp^{e-k}d} \ .\end{equation}

(ii) If $p=2$ and $e\geq3$, then
\begin{eqnarray}\label{e1.7}
{\cal N}(G, n)&=&
 \sum _{d=1}^2\sum_{k=2}^{e} \frac {\varphi (2^{e-k})} {2^{e-1}}\sum\limits_{\lambda_{1}+2\lambda_{2}+\cdots
+n\lambda_{n}=n}\frac{4^{f_2(\lambda_1,\ldots,\lambda_n,\vec{\delta}_e^{(k,d)})}}
{\lambda_1!\cdots\lambda_n!1^{\lambda_1}\cdots n^{\lambda_n}}
\end{eqnarray}
\begin{eqnarray}\label{e1.7'}
&=& \frac{1}{2^{e-1}}\sum\limits_{\lambda_{1}+\cdots
+n\lambda_{n}=n}\frac{4^{e\sum\limits_{t=1}^{n}\lambda_t}+
4^{\sum\limits_{(t,2)=1}\lambda_t+e\sum\limits_{t=1}^{[n/2]}\lambda_{2t}}}
{\lambda_1!\cdots\lambda_n!1^{\lambda_1}\cdots
n^{\lambda_n}}\nonumber\\
&+&\sum\limits_{d=1}^2\sum\limits_{k=2}^{e-1}2^{-k}\sum_{\lambda_{1}+\cdots
+n\lambda_{n}=n}\frac{4^{f_2(\lambda_1,\ldots,\lambda_n,\vec{\delta}_e^{(k,d)})}}
{\lambda_1!\cdots\lambda_n!1^{\lambda_1}\cdots
n^{\lambda_n}},\end{eqnarray} where
\begin{eqnarray}\label{e1.12}
{}&{}
&f_2(\lambda_1,\ldots,\lambda_n,\vec{\delta}_e^{(k,1)})=\sum\limits_{s=0}^{e-k-1}(k+s)\sum\limits_{1\leq
t\leq [n/(2^s)]\atop (t,2)=1}
\lambda_{t2^s}+e\sum\limits_{t=1}^{[n/2^{e-k}]} \lambda_{t2^{e-k}}\
\ ( 2 \le k \le e),  \nonumber\\
{}&{}&
\end{eqnarray}
\begin{eqnarray}\label{e1.13}
f_2(\lambda_1,\ldots,\lambda_n,\vec{\delta}_e^{(k,2)})&=&\sum_{(t,2)=1}\lambda_t+\sum\limits_{s=1}^{e-k-1}(k+s)
\sum\limits_{1\leq t\leq [n/(2^s)]\atop (t,2)=1} \lambda_{t2^s} \nonumber\\
{}&{}& +e\sum\limits_{t=1}^{[n/2^{e-k}]} \lambda_{t2^{e-k}}\ \ \ (
2 \le k \le e-1)
\end{eqnarray}
and
\begin{equation}\label{1.14}
f_2(\lambda_1,\ldots,\lambda_n,\vec{\delta}_e^{(e,2)})=\sum_{(t,2)=1}\lambda_t+
e\sum_{t=1}^{[n/2]}\lambda_{2t}.
\end{equation}

\end{Theorem}
 {\bf Proof.}(i) For $(\phi _i , \sigma )\in
M$, $i \in \mathbb{Z}_{p^e}^*$, by  Lemma\ref {1.3}, there exists a
positive integer  $k$ such that  $i = \alpha ^\nu$, $\nu = p
^{k-1}c$, $(c, p)=1$ and $1 \le k \le e$. By  Lemma\ref {1.4}, there
exists positive integer  $d$ such that  $\vec \delta _e (i) = \vec
\delta _e ^{(k, d)}$ and   $d \mid (p-1)$. First we compute the
number $\mid \! F_{\phi _i, \sigma} \!\mid$ of fixed point set.

$ (1^\circ)$ If $d\nmid r$, then $\delta_1(i)\nmid r$ and $p\nmid
(i^r-1)$, i.e. $(i^r-1,p^e)=1$. By (\ref{e1.5}$'$), $k_{j_1}=0$.

$ (2^\circ)$ If $p^sd|r$ and $p^{s+1}d\nmid r$ ( $0 \le s \le
e-k-1$), then $r=tp^sd$, where $(t,p)=1$. By  Lemma\ref{1.4}, the
order of  $i$ is  $p^sd$ in $\mathbb{Z}_{p^{k+s}}^*$. Consequently,
$i^r-1=cp^{k+s}$. If  $p\! \mid \! c$, i.e. $c=c'p$, then
$i^r-1=c'p^{k+s+1}$.  Since the order of $i$ is $p^{s+1}d$ in
$\mathbb{Z}_{p^{k+s+1}}^*$,  $p^{s+1}d \mid r$. This is a
contradiction. Thus $(c,p)=1$. This implies that the number of
solutions of (\ref{e1.5}$'$) is $p^{k+s}$.

$(3^\circ)$ If $p^{e-k}d|r$, then the number of solutions of
(\ref{e1.5}$'$) is $p^{e}$.

Obviously, (\ref{e1.5}$'$) of all cycles in independent cycle
decomposition of $\sigma$ are independent each other. Notice that
there exactly exist $\sum\limits_{1\leq t\leq [n/(p^sd)]\atop
(t,p)=1} \lambda_{tp^sd}$ cycles, which satisfy $(2^\circ)$ with the
length  $r$, in independent cycle decomposition of $\sigma$; there
exactly exist $\sum\limits_{t=1}^{[n/(p^{e-k}d)]}
\lambda_{tp^{e-k}d}$ cycles, which satisfy $(3^\circ)$ with the
length $r$, in independent cycle decomposition of $\sigma$;
Consequently, there exactly exist
\begin{eqnarray*}
\left(\prod_{s=0}^{e-k-1}p^{(k+s)\sum\limits_{1\leq t\leq
[n/(p^sd)]\atop (t,p)=1}
\lambda_{tp^sd}}\right)p^{e\sum\limits_{t=1}^{[n/p^{e-k}d]}
\lambda_{tp^{e-k}d}}
&=&p^{\sum\limits_{s=0}^{e-k-1}(k+s)\sum\limits_{1\leq t\leq
[n/(p^sd)]\atop (t,p)=1}
\lambda_{tp^sd}+e\sum\limits_{t=1}^{[n/p^{e-k}d]}
\lambda_{tp^{e-k}d}} \\
&=& p^{f_p(\lambda_1,\ldots,\lambda_n,\vec{\delta}_e^{(k,d)})}
\end{eqnarray*} distinct $(k_j)_{j\in J}$ satisfying  (\ref{e1.3}).

Obviously,
\begin{equation}\label{e1.9}
f_p(\lambda_1,\ldots,\lambda_n,\vec{\delta}_e^{(e,d)})=e\sum\limits_{t=1}^{[n/d]}
\lambda_{td}.
\end{equation} If $d>n$, then
\begin{equation}\label{e1.10}
f_p(\lambda_1,\ldots,\lambda_n,\vec{\delta}_e^{(k,d)})=0 \ .
\end{equation}

Similarly, there exactly exist $$
p^{f_p(\lambda_1,\ldots,\lambda_n,\vec{\delta}_e^{(k,d)})}
$$ distinct $(l_j)_{j\in J}$ satisfying (\ref{e1.4}).

Considering the independence  between (\ref{e1.3}) and (\ref{e1.4}),
we have
 \begin{equation}\label{e1.11}
|F_{(\phi_i,\sigma)}|=\left(p^{f_p(\lambda_1,\ldots,\lambda_n,\vec{\delta}_e^{(k,d)})}\right)^2=
p^{2f_p(\lambda_1,\ldots,\lambda_n,\vec{\delta}_e^{(k,d)})}.
\end{equation}

By  Lemma\ref{0.2} and  Lemma\ref{1.4}, there exist
$\varphi(p^{e-k}d) $
$\frac{n!}{\lambda_1!\lambda_2!\cdots\lambda_n!1^{\lambda_1}2^{\lambda_2}\cdots
n^{\lambda_n}}$ elements, whose the number of elements in  fixed
point set is
  $p^{2f_p(\lambda_1,\ldots,\lambda_n,\vec{\delta}_e^{(k,d)})}$, in
 $M$.

 See
\begin{eqnarray*}
{\cal N}(G, n)&=&\frac{1}{|M|}\sum_{(\varphi_i,\sigma)\in
M}|F_{(\varphi_i,\sigma)}|
\end {eqnarray*}
\begin{eqnarray*}
&=&\frac{1}{n!|\mathrm{Aut}G|}\sum_{k=1}^e\sum_{d|p-1}\sum_{\lambda_1+2\lambda_{2}+\cdots+n\lambda_n=n}
\frac{\varphi(p^{e-k}d)n!}{\lambda_1!\cdots\lambda_n!1^{\lambda_1}\cdots
n^{\lambda_n}}p^{2f_p(\lambda_1,\ldots,\lambda_n,\vec{\delta}_e^{(k,d)})}\\
\end{eqnarray*}\begin{eqnarray*}
&=&\frac{1}{p^{e-1}(p-1)}\left(\sum_{d|p-1}\varphi(d)\sum_{\lambda_1+2\lambda_{2}+\cdots+n\lambda_n=n}
\frac{p^{2f_p(\lambda_1,\ldots,\lambda_n,\vec{\delta}_e^{(e,d)})}}{\lambda_1!\cdots\lambda_n!1^{\lambda_1}\cdots
n^{\lambda_n}}\right.\\
&&\left.+\sum_{k=1}^{e-1}p^{e-k-1}(p-1)\sum_{d|p-1}\varphi(d)\sum_{\lambda_1+\cdots+n\lambda_n=n}
\frac{p^{2f_p(\lambda_1,\ldots,\lambda_n,\vec{\delta}_e^{(k,d)})}}{\lambda_1!\cdots\lambda_n!1^{\lambda_1}\cdots
n^{\lambda_n}}\right)\\
&=&\frac{1}{p^{e-1}(p-1)}\sum_{d|p-1}\varphi(d)\sum_{\lambda_1+2\lambda_{2}+\cdots+n\lambda_n=n}
\frac{p^{2f_p(\lambda_1,\ldots,\lambda_n,\vec{\delta}_e^{(e,d)})}}{\lambda_1!\cdots\lambda_n!1^{\lambda_1}\cdots
n^{\lambda_n}}(\hbox{written as (I)})\\
&&+\sum_{k=1}^{e-1}p^{-k}\sum_{d|p-1}\varphi(d)\sum_{\lambda_1+2\lambda_{2}+\cdots+n\lambda_n=n}
\frac{p^{2f_p(\lambda_1,\ldots,\lambda_n,\vec{\delta}_e^{(k,d)})}}{\lambda_1!\cdots\lambda_n!1^{\lambda_1}\cdots
n^{\lambda_n}}(\hbox{written as (II)}).
\end{eqnarray*}
We now compute (I) and (II), respectively. Applying Lemma\ref{0.2},
(\ref{e1.9}) and (\ref{e1.10}), we have
\begin{eqnarray*}
\textrm{(I) }&=&\frac{1}{p^{e-1}(p-1)}\left(\sum_{d|p-1\atop d\leq
n}\varphi(d)\sum_{\lambda_1+2\lambda_{2}+\cdots+n\lambda_n=n}\frac{p^{2e\sum\limits_{t=1}^{[n/d]}\lambda_{td}}}
{\lambda_1!\cdots\lambda_n!1^{\lambda_1}\cdots
n^{\lambda_n}}\right.\\
&&\left.+\sum_{d|p-1\atop d>
n}\varphi(d)\sum_{\lambda_1+2\lambda_{2}+\cdots+n\lambda_n=n}\frac{1}{\lambda_1!\cdots\lambda_n!1^{\lambda_1}\cdots
n^{\lambda_n}}\right)\\
\end{eqnarray*}
\begin{eqnarray*}
&=&\frac{1}{p^{e-1}(p-1)}\left(\sum_{d|p-1\atop d\leq
n}\varphi(d)\sum_{\lambda_1+2\lambda_{2}+\cdots+n\lambda_n=n}\frac{p^{2e\sum\limits_{t=1}^{[n/d]}\lambda_{td}}}
{\lambda_1!\cdots\lambda_n!1^{\lambda_1}\cdots
n^{\lambda_n}}+\sum_{d|p-1\atop d> n}\varphi(d)\right)\\
&=&\frac{1}{p^{e-1}(p-1)}\left(\sum_{d|p-1\atop d\leq
n}\varphi(d)\sum_{\lambda_1+2\lambda_{2}+\cdots+n\lambda_n=n}\frac{p^{2e\sum\limits_{t=1}^{[n/d]}\lambda_{td}}}
{\lambda_1!\cdots\lambda_n!1^{\lambda_1}\cdots
n^{\lambda_n}}+p-1-\sum_{d|p-1\atop d\leq n}\varphi(d)\right)\\
&=&\frac{1}{p^{e-1}(p-1)}\sum\limits_{d|p-1\atop d\leq
n}\left(\sum\limits_{\lambda_{1}+2\lambda_{2}+\cdots
+n\lambda_{n}=n} \varphi
(d)\frac{p^{2e\sum\limits_{t=1}^{[n/d]}\lambda_{td}}}{\lambda_1!\cdots\lambda_n!1^{\lambda_1}\cdots
n^{\lambda_n}}-1\right)+\frac{1}{p^{e-1}}.
\end{eqnarray*}

\begin{eqnarray*}
\textrm{(II)}&=&\sum_{k=1}^{e-1}p^{-k}\left(\sum_{d|p-1\atop d\leq
n}\varphi(d)\sum_{\lambda_1+2\lambda_{2}+\cdots+n\lambda_n=n}
\frac{p^{2f_p(\lambda_1,\ldots,\lambda_n,\vec{\delta}_e^{(k,d)})}}{\lambda_1!\cdots\lambda_n!1^{\lambda_1}\cdots
n^{\lambda_n}}\right.\\
&&\left.+\sum_{d|p-1\atop d>
n}\varphi(d)\sum_{\lambda_1+2\lambda_{2}+\cdots+n\lambda_n=n}
\frac{1}{\lambda_1!\cdots\lambda_n!1^{\lambda_1}\cdots
n^{\lambda_n}}\right)\\
&=&\sum_{k=1}^{e-1}p^{-k}\left(\sum_{d|p-1\atop d\leq
 n}\varphi(d)\sum_{\lambda_1+2\lambda_{2}+\cdots+n\lambda_n=n}
\frac{p^{2f_p(\lambda_1,\ldots,\lambda_n,\vec{\delta}_e^{(k,d)})}}{\lambda_1!\cdots\lambda_n!1^{\lambda_1}\cdots
n^{\lambda_n}}+p-1-\sum_{d|p-1\atop
d\leq  n}\varphi(d)\right)\\
&=&\sum_{k=1}^{e-1}p^{-k}\sum_{d|p-1\atop d\leq
 n}\varphi(d)\left(\sum_{\lambda_1+2\lambda_{2}+\cdots+n\lambda_n=n}
\frac{p^{2f_p(\lambda_1,\ldots,\lambda_n,\vec{\delta}_e^{(k,d)})}}{\lambda_1!\cdots\lambda_n!1^{\lambda_1}\cdots
n^{\lambda_n}}-1\right)+1-\frac{1}{p^{e-1}}.
\end{eqnarray*}

Since ${\cal N}(G, n)=\textrm{I}+\textrm{II}$, (\ref{e1.13}) holds.

(ii)  We first show
\begin{equation} \label {e1.12'}
|F_{(\phi_i,\sigma)}|=\left(2^{f_2(\lambda_1,\ldots,\lambda_n,\vec{\delta}_e^{(k,d)})}\right)^2=
2^{2f_2(\lambda_1,\ldots,\lambda_n,\vec{\delta}_e^{(k,d)})}.
\end{equation}

Case 1:   $i = \beta  ^\nu$ and  $\nu = 2 ^{k-2}c$ with  $(c, p)=1$
and $k
>1$. By  Lemma\ref {1.4}, $d=1$ and   $\vec \delta _e (i) =
\vec \delta _e ^{(k, 1)}$.

$(1^\circ)$ If  $2 ^s \mid r$ and  $2^{s+1}\nmid r$ with $0\leq
s\leq e-k-1)$, then the number of solutions of (\ref{e1.5}$'$) is
$2^{k+s}$ by means of the method similar to proof of (i).

$(2^\circ)$ If $2^{e-k}|r$, then the number of solutions of
(\ref{e1.5}$'$) is $2^e$.

There exists $\sum\limits_{1\leq t\leq [n/(2^s)]\atop (t,2)=1}
\lambda_{t2^s}$ cycles, which satisfy $(1^\circ)$ with the length
$r$, in independent cycle decomposition of $\sigma$; there exist
$\sum\limits_{t=1}^{[n/2^{e-k}]} \lambda_{t2^{e-k}}$ cycles, which
satisfy $(2^\circ)$ with the length $r$, in independent cycle
decomposition of $\sigma$. This implies (\ref{e1.12'}).

{\bf Remark.} It is clear
\begin{equation}\label{e2.11}
f_2(\lambda_1,\ldots,\lambda_n,\vec{\delta}_e^{(e,1)})
=e\sum\limits_{t=1}^{n} \lambda_{t}.
\end{equation}

Case 2:   $i = - \beta ^\nu$ and  $\nu = 2 ^{k-2}c$ with  $(c, p)=1$
and  $e \ge k
>1$. By  Lemma\ref {1.4}, $d=2$ and   $\vec \delta _e (i) =
\vec \delta _e ^{(k, 2)}$.

$(1^\circ)$ If $2\nmid r$, i.e. $r$ is  odd. It is clear $i^r-1=2c$.
Since  $\delta _4 (i) =2$,  $c$ is odd. By (\ref{e1.5}$'$),
$k_{j_1}=c'2^{e-1}$. Thus the number of solutions of (\ref{e1.5}$'$)
is 2.

$(2^\circ)$ If $2^s|r$ and $2^{s+1}\nmid r\ (1\leq s\leq e-k-1)$,
then $r=t2^s$ and $(t,2)=1$;  $i^r-1=c2^{s+k}$ with $(c,2)=1$;
$k_{j_1}=c'2^{e-s-k}$. Thus the number of solutions of
(\ref{e1.5}$'$) is $2^{s+k}$.

$(3^\circ)$ If $2^{e-k}|r$, then  the number of solutions of
(\ref{e1.5}$'$) is $2^{e}$.

There exists $\sum_{(t,2)=1}\lambda_t$
 cycles, which satisfy $(1^\circ)$ with the length
$r$, in independent cycle decomposition of $\sigma$; there exist
$\sum\limits_{1\leq t\leq [n/(2^s)]\atop (t,2)=1} \lambda_{t2^s}$
cycles, which satisfy $(2^\circ)$ with the length $r$, in
independent cycle decomposition of $\sigma$; there exist
$\sum\limits_{t=1}^{[n/2^{e-k}]} \lambda_{t2^{e-k}}$ cycles, which
satisfy $(3^\circ)$ with the length $r$, in independent cycle
decomposition of $\sigma$. This implies (\ref{e1.12'}).

Now we show (ii). Applying  Lemma\ref{1.3}(ii), Lemma\ref{1.4} and
Burnside's Lemma, we have
\begin{eqnarray*}
{\cal N}(G, n)&=&\frac{1}{|M|}\sum_{(\varphi_i,\sigma)\in
M}|F_{(\varphi_i,\sigma)}|\\
&=&\frac{1}{n!|\mathrm{Aut}G|}\left( \sum
_{d=1}^2\sum_{k=2}^e\sum_{\lambda_1+2\lambda_{2}+\cdots+n\lambda_n=n}
\varphi(2^{e-k})\frac{n!}{\lambda_1!\cdots\lambda_n!1^{\lambda_1}\cdots
n^{\lambda_n}}4^{f_2(\lambda_1,\ldots,\lambda_n,\vec{\delta}_e^{(k,d)})}
\right. \\
\end{eqnarray*}
\begin{eqnarray*}
&=&\frac{1}{2^{e-1}}\left(\sum_{d=1}^2\sum_{k=2}^{e-1}2^{e-k-1}
\sum_{\lambda_1+\cdots+n\lambda_n=n}
\frac{4^{f_2(\lambda_1,\ldots,\lambda_n,\vec{\delta}_e^{(k,d)})}}
{\lambda_1!\cdots\lambda_n!1^{\lambda_1}\cdots
n^{\lambda_n}}\right.\\
 &&\left.+\sum_{\lambda_1+\cdots+n\lambda_n=n}\frac{4^{e\sum\limits_{t=1}^{n}
\lambda_{t}}} {\lambda_1!\cdots\lambda_n!1^{\lambda_1}\cdots
n^{\lambda_n}}+ \sum_{\lambda_1+\cdots+n\lambda_n=n}
\frac{4^{\sum\limits_{(t,2)=1}\lambda_t+
e\sum\limits_{t=1}^{[n/2]}\lambda_{2t}}}{\lambda_1!\cdots\lambda_n!1^{\lambda_1}\cdots
n^{\lambda_n}}\right)\\
\end{eqnarray*}
\begin{eqnarray*}
&=&\sum_{d=1}^2\sum_{k=2}^{e-1}2^{-k}\sum\limits_{\lambda_{1}+\cdots
+n\lambda_{n}=n}\frac{4^{f_2(\lambda_1,\ldots,\lambda_n,\vec{\delta}_e^{(k,d)})}}
{\lambda_1!\cdots\lambda_n!1^{\lambda_1}\cdots n^{\lambda_n}}\\
&&+\frac{1}{2^{e-1}}\sum\limits_{\lambda_{1}+\cdots
+n\lambda_{n}=n}\frac{4^{e\sum\limits_{t=1}^{n}\lambda_t}
+4^{\sum\limits_{(t,2)=1}\lambda_t+e\sum\limits_{t=1}^{[n/2]}\lambda_{2t}}}
{\lambda_1!\cdots\lambda_n!1^{\lambda_1}\cdots n^{\lambda_n}}.\ \
 \ \ \Box  \nonumber\\
\end{eqnarray*}

\begin{Corollary}\label{1.7}
Let $G\cong C_{p^e}$.

(i) If $p$ is an odd prime number or $e\leq2$, then
$$
{\cal N}(G, 1)=p^e+2(p^{e-1}+p^{e-2}+\cdots+p+1);
$$
(ii)  If $p=2$ and $e\geq3$, then
$${\cal N}(G, 1)=2^{e+1}+2^e-2.$$
\end{Corollary}

 {\bf Proof.} (i)
Notice  $f(\lambda_1,\vec{\delta}_e^{(k,1)})=k (1\leq k\leq e-1)$
and  $f(\lambda_1,\vec{\delta}_e^{(e,1)})=e$  when  $n=1$. By
(\ref{e1.6}), we have
\begin{eqnarray*}
{\cal N}(G,1)&=&1+\frac{1}{p^{e-1}(p-1)}\left(p^{2e}-1\right)+\sum_{k=1}^{e-1}p^{-k}\left(p^{2k}-1\right)\\
&=&1+\frac{1}{p^{e-1}(p-1)}\left(p^{2e}-1\right)+\sum_{k=1}^{e-1}p^k-\sum_{k=1}^{e-1}p^{-k}\\
&=&p^e+2(p^{e-1}+p^{e-2}+\cdots+p+1).
\end{eqnarray*}

(ii) By (\ref{e1.12}) and  (\ref{e1.13}),
$f(\lambda_1,\vec{\delta}_2^{(k,1)})=k (2\leq k\leq e-1)$ and
$f(\lambda_1,\vec{\delta}_2^{(k,2)})=1 (2\leq k\leq e-1)$ when
$n=1$. Applying (\ref{e1.7'}), we have
\begin{eqnarray*}
{\cal
N}(G,1)&=&\frac{1}{2^{e-1}}(4^e+4)+\sum_{k=2}^{e-1}(2^{-k}4^k)+\sum_{k=2}^{e-1}(2^{-k}4)\\
&=&2^{e+1}+2^e-2. \Box
\end{eqnarray*}

\begin{Corollary}\label{1.77}
 Assume that  $G=C_{p^e}$ and  $p$ is a  prime number with positive integer $e$.

(i) If  $p$ is an odd prime number, then $$ {\cal N}(G, 2
)=1+\frac{1}{2}\frac{p^{3e}-p^3}{p^3-1}+\frac{1}{p-1}\left(\frac{1}{2}p^{3e+1}
+p^{e+1}+p^e-p-\frac{3}{2}\right);
$$

(ii) ${\cal N}(C_4, 2 )=76, \  {\cal N}(C_2, 2 )=10$;

(iii) If $p=2$ and $e\geq3$, then
$${\cal N}(G, 2)=\frac{15}{14}2^{3e}+3\times 2^{e+1}-\frac{116}{7}.$$

\end{Corollary}

{\bf Proof.} It is clear that $(\lambda_1,\lambda_2)=(0,1)$ or (2,0)
when $n=2$.

(i) If $p$ is an odd prime number, by (\ref{e1.8}), we have $d =
1, 2$,  and $f_p(\lambda_1,\lambda_2,\vec{\delta}_e^{(k,1)})$
$=k(\lambda_1+\lambda_2)$ and $
f_p(\lambda_1,\lambda_2,\vec{\delta}_e^{(k,2)})$ $=k\lambda_2
(1\leq k\leq e)$.  Using (\ref{e1.6'}), we have
\begin{eqnarray*}
{\cal
N}(G,2)&=&1+\frac{1}{p^{e-1}(p-1)}\left(\sum_{\lambda_1,\lambda_2}
\frac{p^{2e(\lambda_1+\lambda_2)}}
{\lambda_1!\lambda_2!1^{\lambda_1}
2^{\lambda_2}}-1+\sum_{\lambda_1,\lambda_2} \frac{p^{2e\lambda_2}}
{\lambda_1!\lambda_2!1^{\lambda_1} 2^{\lambda_2}}-1\right)+\\
&&\sum_{k=1}^{e-1}p^{-k}\left(\sum_{\lambda_1,\lambda_2}
\frac{p^{2k(\lambda_1+\lambda_2)}}
{\lambda_1!\lambda_2!1^{\lambda_1}
2^{\lambda_2}}-1+\sum_{\lambda_1,\lambda_2} \frac{p^{2k\lambda_2}}
{\lambda_1!\lambda_2!1^{\lambda_1} 2^{\lambda_2}}-1\right)\\
 &=&1+\frac{1}{p^{e-1}(p-1)}\left(\frac{1}{2}p^{2e}+\frac{1}{2}p^{4e}-1+\frac{1}{2}p^{2e}
 +\frac{1}{2}-1\right)+\\
&&\sum_{k=1}^{e-1}p^{-k}\left(\frac{1}{2}p^{2k}+\frac{1}{2}p^{4k}-1+\frac{1}{2}p^{2k}
 +\frac{1}{2}-1\right)\\
 &=&1+\frac{1}{2}\frac{p^{3e}-p^3}{p^3-1}+\frac{1}{p-1}\left(\frac{1}{2}p^{3e+1}
+p^{e+1}+p^e-p-\frac{3}{2}\right).
\end{eqnarray*}

(ii) If $p=2$ and $e=2$, then $d=1$ and $k=1$. By (\ref{e1.6'}), we
have
\begin{eqnarray*}
{\cal N}(C_4, 2 )&=&1+\frac{1}{2}\left(\sum_{\lambda_1,\lambda_2}
\frac{4^{e(\lambda_1+\lambda_2)}} {\lambda_1!\lambda_2!1^{\lambda_1}
2^{\lambda_2}}-1\right)+\frac{1}{2}\left(\sum_{\lambda_1,\lambda_2}
\frac{4^{(\lambda_1+2\lambda_2)}}
{\lambda_1!\lambda_2!1^{\lambda_1} 2^{\lambda_2}}-1\right)\\
&=&76.
\end{eqnarray*}

If $p=2$ and $e=1$, then $d=1$. By (\ref{e1.6}), we have
\begin{eqnarray*}
{\cal N}(C_2, 2)&=&1+\left(\sum_{\lambda_1,\lambda_2}
\frac{4^{\lambda_1+\lambda_2}}{\lambda_1!\lambda_2!1^{\lambda_1}2^{\lambda_2}}-1\right)\\
&=&10.
\end{eqnarray*}

(iii) It follows from (\ref{e1.13}) and (\ref{1.14}) that
$f_2(\lambda_1,\lambda_2,\vec{\delta}_e^{(k,1)})=k\lambda_1+(k+1)\lambda_2$
and $
f(\lambda_1,\lambda_2,\vec{\delta}_e^{(k,2)})=\lambda_1+(k+1)\lambda_2
(2\leq k\leq e-1)$. Applying (\ref{e1.7'}), we have
\begin{eqnarray*}
{\cal N}(G, 2)&=&\frac{1}{2^{e-1}}\sum_{\lambda_1,\lambda_2}
\frac{4^{e(\lambda_1+\lambda_2)}+4^{\lambda_1+e\lambda_2}}
{\lambda_1!\lambda_2!1^{\lambda_1}2^{\lambda_2}}+\sum_{k=2}^{e-1}2^{-k}
\sum_{\lambda_1,\lambda_2} \frac{4^{k\lambda_1+(k+1)\lambda_2}}
{\lambda_1!\lambda_2!1^{\lambda_1}2^{\lambda_2}} \\
&&+\sum_{k=2}^{e-1}2^{-k} \sum_{\lambda_1,\lambda_2}
\frac{4^{\lambda_1+k\lambda_2}}
{\lambda_1!\lambda_2!1^{\lambda_1}2^{\lambda_2}}\\
&=&\frac{1}{2^{e-1}}\left(\frac{4^{2e}+4^2}{2}+\frac{4^e+4^e}{2}\right)
+\sum_{k=2}^{e-1}2^{-k}\left(\frac{4^{2k}+4^{k+1}}{2}\right)\\
&&+\sum_{k=2}^{e-1}2^{-k}\left(\frac{4^2+4^{k+1}}{2}\right)\\
&=&\frac{15}{14}2^{3e}+3\times 2^{e+1}-\frac{116}{7}. \Box
\end{eqnarray*}

\begin{Corollary}\label{1.8}
 If $G\cong C_p$ and  $p$ is a prime number, then
$${\cal N}(G,
n)=1+\frac{1}{p-1}\sum\limits_{d|p-1\atop d\leq n}
\varphi(d)\left(\sum\limits_{\lambda_{1}+2\lambda_{2}+\cdots
+n\lambda_{n}=n}\frac{p^{2\sum\limits_{t=1}^{[n/d]}
\lambda_{td}}}{\lambda_1!\cdots\lambda_n!1^{\lambda_1}\cdots
n^{\lambda_n}}-1\right).
$$
\end{Corollary}

{\bf Proof.} If follows from (\ref{e1.6'}).\ $\Box$

\section{Finite Cycle Groups}\label {s2}

In this section we give the formula computing the number of
isomorphic classes of element systems over finite cycle groups.

\begin{Lemma} \label{2.1}
 Assume $G=C_m$ and $m=p_1^{e_1}p_2^{e_2}\cdots
p_s^{e_s}$, where $p_1,\ldots,p_s$ are mutually different prime
numbers, then $G=G_1\oplus\cdots\oplus
 G_s$ with $ G_i \cong  C_{m_i}$ and  $m_i = p_i^{e_i}$. Furthermore,

(i) Aut$G\cong$ Aut$G_1\oplus\cdots \oplus$Aut$G_s$;

(ii) There exists a bijective map $\psi:\widehat{G}\rightarrow
\widehat{G_1}\oplus\cdots\oplus\widehat{G_s}$.

\end{Lemma}

{\bf Proof.} (i) It follows from  \cite [Theorem 1.11.10] {Zh82}.

(ii)  Define map $\psi $:  $\widehat{G}\rightarrow
\widehat{G_1}\oplus\cdots\oplus\widehat{G_s}$ by sending $\chi \in
\widehat{G}$    to
$$\psi(\chi)=(\chi_{1},\ldots,\chi_{s}),$$
where $\chi_{i}$ is the restriction of $\chi$ on $G_i$. It is clear
that $\psi$ is injective.

Conversely, for any $(\chi_1,\ldots,\chi_s)\in
\widehat{G_1}\oplus\cdots\oplus\widehat{G_s}$,
 Define  $\chi \in \hat G$ such that
$ \chi(g)=\chi_1(g_1)\cdots\chi_s(g_s),$ for any
$g=(g_1,\ldots,g_s)\in G,g_i\in G_i.$ Thus $\psi$ is surjective.
$\Box$

\begin{Theorem}\label{2.2}
 Assume that finite cycle group $G\cong C_{p_1^{e_1}}\oplus
C_{p_2^{e_2}}\oplus \cdots \oplus C_{p_s^{e_s}}$ and $p_i$ is a
prime number with positive integer  $e_i$ and  $p_1< p_2< \cdots
<p_s.$

(i) If all of  $p_1,\cdots,p_s$ are odd prime numbers or $p_1=2$
with $e_1\leq 2$, then
\begin{equation}\label{e7.6}
{\cal N}(G, n)=\sum_{\lambda_1 + \cdots + n \lambda _n
=n}\left(\frac{1 }{\lambda_1!\cdots\lambda_n!1^{\lambda_1}\cdots
n^{\lambda_n}}\prod_{i=1}^{s
}\frac{\sum\limits_{k_i=1}^{e_i}\sum\limits_{d_i|p_i-1}
\varphi(p_i^{e_i-k_i}d_i)p_i^{2f_{p_i}(\lambda_1, \ldots,\lambda_n,
\vec \delta_{e_i}^{(k_i,d_i)})} }{p_i^{e_i-1}(p_i-1)}\right) .
\end{equation}

(ii) If $p_1=2$ and  $e_1 \ge 3$,  then
\begin{eqnarray*}
{\cal N}(G, n)&=& \sum_{\lambda_1 + \cdots + n \lambda _n
=n}\left(\frac{1 }{\lambda_1!\cdots\lambda_n!1^{\lambda_1}\cdots
n^{\lambda_n}}\frac{\sum\limits_{k_1=2}^{e_1}\sum\limits_{d_1=1}^{2}
\varphi(2^{e_1-k_1})4^{f_{2}(\lambda_1,
\ldots,\lambda_n, \vec \delta_{e_1}^{(k_1,d_1)})}}{2^{e_1-1}}\right.\\
&&\left.\prod_{i=2}^{s
}\frac{\sum\limits_{k_i=1}^{e_i}\sum\limits_{d_i|p_i-1}
\varphi(p_i^{e_i-k_i}d_i)p_i^{2f_{p_i}(\lambda_1,
\ldots,\lambda_n,\vec \delta_{e_i}^{(k_i,d_i)})}
}{p_i^{e_i-1}(p_i-1)}\right) .
\end{eqnarray*}

\end{Theorem}
{\bf Proof.} Given $(\phi,\sigma)\in M$, we first compute
$|F_{(\phi,\sigma)}|$. According to  Lemma\ref{2.1}, we can assume
$\phi=(\phi_1,\ldots,\phi_s)$ with $\phi_i\in \mathrm{Aut}G_i$. It
is clear
\begin{equation}\label{e7.1}
\phi(g)=(\phi(g_1),\ldots,\phi(g_s))
\end{equation}
and
\begin{equation}\label{e7.2}
\chi\phi^{-1}=(\chi_1\phi_1^{-1},\ldots,\chi_s\phi_s^{-1})
\end{equation}
for any $g=(g_1,\ldots,g_s)\in G$, $\chi=(\chi_1,\ldots,\chi_s)\in
\widehat{G}$.

If $\left(\begin{array}{c} j\\
g_{j}\\
\chi_{j} \end{array}\right)\in F_{(\phi,\sigma)}$, applying
(\ref{e0.1}), (\ref{e7.1}) and (\ref{e7.2}), we have
\begin{eqnarray}\label{e7.3}
 (\phi,\sigma)\circ \left(\begin{array}{c} j\\
g_{j}\\
\chi_{j} \end{array}\right)
&=&\left(\begin{array}{c} \sigma(j)\\
\phi(g_{j})\\
\chi_{j}\phi^{-1} \end{array}\right) \nonumber
=\left(\begin{array}{c} \sigma(j)\\
(\phi(g_{j1}),\ldots,\phi(g_{js}))\\
(\chi_{j1}\phi_1^{-1},\ldots,\chi_{js}\phi_s^{-1})
\end{array}\right)\\ &=&\left(\begin{array}{c}j\\
(g_{j1},\ldots,g_{js})\\
(\chi_{j1},\ldots,\chi_{js}) \end{array}\right),
\end{eqnarray}
where  $g_j=(g_{j1},\ldots,g_{js}),
\chi_j=(\chi_{j1},\ldots,\chi_{js})$. (\ref{e7.3})is equivalent to
\begin{eqnarray*}
\left(\begin{array}{c} \sigma(j)\\
\phi_i(g_{ji})\\
(\chi_{ji}\phi_i^{-1})
\end{array}\right)=\left(\begin{array}{c}j\\
g_{ji}\\
\chi_{ji}\end{array}\right),\ 1\leq i\leq s.
\end{eqnarray*}
This implies  $\left(\begin{array}{c}j\\
g_{ji}\\
\chi_{ji}\end{array}\right)\in F_{(\phi_i,\sigma)}$, where
$F_{(\phi_i,\sigma)}$ denotes the fixed point set of
$(\phi_i,\sigma)$ in $\Omega(G_i,n)$. Consequently,
\begin{equation}\label{e7.4}
|F_{(\phi,\sigma)}|=\prod_{j=1}^{s}|F_{(\phi_{i},\sigma)}|.
\end{equation}

(i) By  (\ref {e1.11}) we have
\begin{eqnarray*} {\cal N}(G,
n)&=&\frac{1}{|M|}\sum_{(\phi,\sigma)\in
M}|F_{(\phi,\sigma)}|\\
&=&\frac{1}{n!|\mathrm{Aut}G|}\sum_{\lambda_1 + \cdots + n \lambda
_n
=n}\sum_{k_1=1}^{e_1}\sum_{d_1|p_1-1}\cdots\sum_{k_s=1}^{e_s}\sum_{d_s|p_s-1}
\left(\frac{n! }{\lambda_1!\cdots\lambda_n!1^{\lambda_1}\cdots
n^{\lambda_n}}\right.\\
&&\left.\prod_{i=1}^{s }\varphi(p_i^{e_i-k_i}d_i)
p_i^{2f_{p_i}(\lambda_1,
\ldots,\lambda_n, \vec \delta_{e_i}^{(k_i,d_i)})}\right)\\
&=&\sum_{\lambda_1 + \cdots + n \lambda _n =n}\left(\frac{1
}{\lambda_1!\cdots\lambda_n!1^{\lambda_1}\cdots
n^{\lambda_n}}\prod_{i=1}^{s
}\frac{\sum\limits_{k_i=1}^{e_i}\sum\limits_{d_i|p_i-1}
\varphi(p_i^{e_i-k_i}d_i)p_i^{2f_{p_i}(\lambda_1,
\ldots,\lambda_n,\vec \delta_{e_i}^{(k_i,d_i)})}
}{p_i^{e_i-1}(p_i-1)}\right).
\end{eqnarray*}

(ii) By  (\ref {e1.11}) and \ref {e1.12'}), we have
\begin{eqnarray*}
{\cal N}(G, n) &=&\frac{1}{n!|\mathrm{Aut}G|}\sum_{\lambda_1 +
\cdots + n \lambda _n
=n}\sum_{k_1=2}^{e_1}\sum_{d_1=1}^2\sum_{k_2=1}^{e_1}\sum_{d_2|p_2-1}\cdots\sum_{k_s=1}^{e_s}\sum_{d_s|p_s-1}
\left(\frac{n! }{\lambda_1!\cdots\lambda_n!1^{\lambda_1}\cdots
n^{\lambda_n}}\right.\\
&&\left.\varphi(2^{e_1-k_1})p_1^{2f_{p_1}(\lambda_1,
\ldots,\lambda_n,\vec \delta_{e_i}^{(k_1,d_1)})}\prod_{i=2}^{s
}\varphi(p_i^{e_i-k_i}d_i)p_i^{2f_{p_i}(\lambda_1,
\ldots,\lambda_n,\vec \delta_{e_i}^{(k_i,d_i)})}\right). \ \ \Box
\end{eqnarray*}

\begin{Corollary} \label {2.3}
 If $G\cong C_{p_1}\oplus\cdots\oplus C_{p_s}$ and  $p_1,
\ldots,p_s$ are mutually different prime numbers, then
$$\mathcal{N}(G,1)=\prod_{i=1}^s(p_i+2).$$
\end{Corollary}
 {\bf Proof.}  Obviously $k_i=e_i=1$,  $ f_{p_i}(\lambda_1,\vec \delta_{e_i}^{(1,1)})=1$ and
 $
f_{p_i}(\lambda_1, \vec \delta_{e_i}^{(1,d_i)})=0(d_i>1)$. Using
(\ref{e7.6}), we have
\begin{eqnarray*}
{\cal N}(G, 1)
&=&\prod_{i=1}^s\frac{\sum\limits_{d_i|p_i-1}\varphi(d_i)p_i^{2f_{p_i}(\lambda_1,
,\vec \delta_{e_i}^{(1,d_i)})}}{p_i-1}\\
&=&\prod_{i=1}^s\frac{p_i^2+\sum\limits_{d_i|p_i-1\atop
d_i>1}\varphi(d_i)}{p_i-1}\\
&=&\prod_{i=1}^s\frac{p_i^2+p_i-2}{p_i-1}\\
&=&\prod_{i=1}^s(p_i+2).\Box
\end{eqnarray*}

\section{Primary Commutative Groups}\label {s3}

In this section we give the formula computing the number of
isomorphic classes of element systems of primary commutative groups.

If $G\cong \overbrace{C_p\times\cdots\times C_p}^{s}$ and
 $p$ is a prime number, then $G$ is called a primary commutative group.
\begin{Lemma}\label{3.1} (See \cite [ Theorem 5.4.2]{Zh82})
 If $G$ is a  primary commutative group with order $p^s$ and positive integer $s$, then

(i) Aut$G\cong GL(s,\mathbb{Z}_p);$

(ii)
$|$Aut$G|$=$p^{\frac{1}{2}s(s-1)}\prod\limits_{i=1}^{s}(p^i-1).$
\\ Here $GL(s,\mathbb{Z}_p)$ is a general linear group over field $\mathbb{Z}_p$.
Indeed, $GL(s,\mathbb{Z}_p)$ also can be viewed as $s \times
S$-matrix ring over field $\mathbb{Z}_p$.
\end{Lemma}

Let  $g_1,\cdots,g_s$ be a basis of $G$, i.e. $G\cong \langle
g_1\rangle\times\cdots\times \langle g_s \rangle$ with $o(g_i)=p,$
where  $o(g_i)$ denotes the order of  $g_i$. For any $ g\in G$,
$g=g_1^{k_1}\cdots g_s^{k_s},$ let $\vec{k}:=\left (
\begin {array} {c}  k_1
\\
\vdots\\
k_s\end {array}  \right ) denote \  $g$.$

By  Lemma\ref {3.1},for any $ \phi \in$Aut$G$, there exists a matrix
$A=(a_{ij})_{s\times s}\in GL(s,\mathbb{Z}_p)$ such that
\begin{align*}
\phi:G &\rightarrow G
\\\vec{k} &\mapsto A \vec{k}.
\end{align*} Let $A$ denote $\phi$.

$\hat{G}$ consists of the following maps:
\begin{align*}\chi:G &\rightarrow F
\\g_i &\mapsto \omega^{l_i}\qquad i=1,2,\cdots,s; \ 0\leq l_i \leq p-1.
\end{align*}\\ Here $\omega$ is a primitive $p$th root of 1. Let $\vec{l}=(l_1,\cdots,l_s)^T$ denote the character
$\chi$ above, where  $T$ denotes  transposition.

Next we compute the fixed point set $F_{(A,\sigma)}$ of
$(A,\sigma)$.
If $\left(\begin{array}{c} j\\
\vec{k}_j\\
\vec{l}_j\end{array}\right)\in F_{(A,\sigma)}$, then
$$\left(\begin{array}{c} \sigma(j)\\
A\vec{k}_j\\
(A^{-1})^{\mbox{\tiny T}}\vec{l}_j \end{array}\right)=\left(\begin{array}{c} j\\
\vec{k}_j\\
\vec{l}_j\end{array}\right)=\left(\begin{array}{c} \sigma(j)\\
\vec{k}_{\sigma(j)}\\
\vec{l}_{\sigma(j)}\end{array}\right).$$  That is,
\begin {eqnarray} \label {e3.1}A\vec{k}_j=\vec{k}_{\sigma(j)} \end  {eqnarray}
and \begin {eqnarray} \label {e3.1'}
 (A^{-1})^{\mbox{\tiny T}}\vec{l}_j=\vec{l}_{\sigma(j)},\qquad
j=1,\cdots,n, \end {eqnarray} where $\vec{k_j},  \vec{l_j}\in (
\mathbb{Z}_p )^s $.

Assume that $\tau=(j_1,\cdots,j_r)$ is a cycle in the independent
cycle decomposition  of $\sigma$ with the length $r$. We have
\begin {eqnarray} \label {e3.2}
\begin{array} {llll}
\vec k_{j_2}&= A  \vec k_{j_1},& \vec l_{j_2}&=B \vec l_{j_1},\\
&\cdots & \cdots & \\
\vec k_{j_r}&=A^{r-1}\vec k_{j_1}, &\vec l_{j_r}&=B^{r-1}\vec l_{j_1},\\
\vec k_{j_1}&=A^r\vec k_{j_1},& \vec l_{j_1}&=B^r\vec l_{j_1},
\end{array}
\end  {eqnarray}  where $B=(A^{-1})^{\mbox{\tiny
T}}$. Thus we only need choice $\vec k_{j_1}$ and $\vec l_{j_1}$
such that $\vec k_{j_1}$ and $\vec l_{j_1}$ become the solutions of
$(A^r-I)X=0$ and $(B^r-I)X=0$, respectively. Applying (\ref {e3.2}),
we compute $\vec k_{j_2},\cdots,\vec k_{j_r},\vec
l_{j_2},\cdots,\vec l_{j_r}$. This implies that we find a  solution
 $ \vec k_j $ ' s , $ \vec l_j{} $'s  of (\ref {e3.1}) for $j = j_1, j_2, \cdots , j_r$. Conversely,  it is
true, too.

Considering
\begin{eqnarray*}
B^r-I&=&((A^r)^{\mbox{\tiny T}})^{-1}-I\\
&=&((A^r)^{\mbox{\tiny T}})^{-1}-((A^r)^{\mbox{\tiny
T}})^{-1}(A^r)^{\mbox{\tiny T}}\\
&=&((A^r)^{\mbox{\tiny T}})^{-1}(I-(A^r)^{\mbox{\tiny T}}),
\end{eqnarray*}
we have rank$(A^r-I)=$rank$(B^r-I)$, i.e. the numbers of solutions
of
 of $(A^r-I)X=0$ and $(B^r-I)X=0$ in $ (\mathbb{Z}_p)^s$ are the same, written as $\nu
_{A, r}$. Thus  \begin {eqnarray}\label {e3.3} \nu _{A, r} = p ^{ s
-rank (A^r-I)}.
\end {eqnarray}

Assume the type of $\sigma$ is $1^{\lambda_1}2^{\lambda_2}\cdots
n^{\lambda_n}$. It is clear
\begin {eqnarray} \label {e3.4}
|F_{(A,\sigma)}|=\prod_{r=1}^{n}\nu _{A, r}^{2\lambda_r}.\end
{eqnarray}

\begin{Theorem}\label{3.2}

If $G$ is a primary commutative group with order $p^s$, then
$$
{\cal N}(G,
n)=\frac{1}{p^{\frac{1}{2}s(s-1)}\prod\limits_{i=1}^{s}(p^i-1)}\sum\limits_{A\in
GL(s,\mathbb{Z}_p)}\ \sum\limits_{\lambda_1+2\lambda_{2}+\cdots
+n\lambda_n=n}\ \prod\limits_{i=1}^{n}\frac{p ^{2( s -rank
(A^i-I)){\lambda_i}}}{\lambda_i!i^{\lambda_i}}.
$$
\end{Theorem}

{\bf Proof.} Applying  Lemma \ref {3.1} and  (\ref {e3.4}), we
have
\begin{eqnarray*}
{\cal N}(G, n)&=&\frac{\sum\limits_{(A,\sigma)\in M}|F_{(A,\sigma)}|}{|M |}\nonumber\\
&=&\frac{\sum\limits_{A\in GL(s,\mathbb{Z}_p)}\
\sum\limits_{\lambda_1+2\lambda_{2}+\cdots
+n\lambda_n=n}\frac{n!}{\lambda_1!\cdots\lambda_n!1^{\lambda_1}\cdots
n^{\lambda_n}}\ \prod\limits_{r=1}^{n}\nu _{A, r}^{2\lambda_r}}{p^{\frac{1}{2}s(s-1)}\prod\limits_{i=1}^s(p^i-1)n!}\nonumber\\
&=&\frac{1}{p^{\frac{1}{2}s(s-1)}\prod\limits_{i=1}^{s}(p^i-1)}\sum\limits_{A\in
GL(s,\mathbb{Z}_p)}\ \sum\limits_{\lambda_1+2\lambda_{2}+\cdots
+n\lambda_n=n}\ \prod\limits_{i=1}^{n}\frac{\nu _{A,
i}^{2\lambda_i}}{\lambda_i!i^{\lambda_i}}.\Box
\end{eqnarray*}

\section{Finite Commutative Groups}\label {s4}
In this section we give the formula computing the number of
isomorphic classes of element systems over finite commutative
groups.

Let  $G$ be an additive group with order $m$ and $h_1,\cdots,h_s$
be a basis of $G$,i.e. $G= <h_1> \times <h_2> \times \cdots \times
<h_s>$, where  the order of $h_i$ is $m_i =p_i^{e_i}$ and $p_i$ is
a prime number with positive integer $e_i$.

%也就是说, $G\cong Z_{m_1}\times Z_{m_2}\times \cdots \times
%Z_{m_s}$. 不妨设 $G= Z_{m_1}\times Z_{m_2}\times \cdots \times
%Z_{m_s}$. 让 $h_i$表示 $Z_{m_i}$ 中的 $\bar 1$,$\bar E$ 表示
%对角矩阵 diag $ (h_1, h_2 , \cdots , h_s) $ $: =\left ( \begin
%{array} {ccccc} h_1 & 0 & \cdots & 0 \\
%0 & h_2 &  \cdots & 0 \\
%\cdots  &\cdots & \cdots & \cdots  \\
%0& 0 & \cdots & h_s \\
%\end {array}\right )$.
%那么, $\forall g\in G, $ $g=h_1{k_1}+ \cdots + h_s{k_s}$ = $(h_1,
%h_2, \cdots, h_s)\left (
%\begin {array} {c}  k_1
%\\
%\vdots\\
%k_s\end {array}  \right ),$ 其中, $k_i \in Z_{m_i}$, \ $i =1, 2,
%\cdots, s.$ 将$g$记为$\vec{k}:=\left (
%\begin {array} {c}  k_1
%\\
%\vdots\\
%k_m\end {array}  \right ).$

Let $M_{s\times s}(\mathbb{Z})$ denote the set of all $s\times s$-
matrices  over $\mathbb{Z}$. For any  $ \phi \in$Aut$G$, there
exists $A=(a_{ij})_{s\times s}\in M_{s\times s} (\mathbb{Z})$,such
that  $\phi(h_1, h_2, \cdots, h_s) = (h_1, h_2, \cdots, h_s) A.$
Similarly, there exists $B=(b_{ij})_{s\times s}\in M_{s\times s}
(\mathbb{Z})$ such that
 $\phi ^{-1}(h_1, h_2, \cdots, h_s)$ $ = $ $(h_1, h_2, \cdots,$ $
h_s) B^{T}$.

 Let $\omega$ be a primitive $m$th root of 1. For   $\vec k  \in  \left (
\begin {array} {l} \mathbb{Z}_{m_1}
\\
\vdots\\
\mathbb{Z}_{m_s} \end {array}  \right)$ and  $\vec l \in  \left (
\begin {array} {l} \mathbb{Z}_{m}
\\
\vdots\\
\mathbb{Z}_{m} \end {array}  \right)$,   define   $g_{\vec k}
:=h_1{k_1}+ \cdots + h_s{k_s} \in G$ and  $\chi_ {\vec l } \in \hat
G$ such that $\chi_ {\vec l } (h_1, h_2, \cdots, h_s) := (\omega
^{l_1},\omega ^{l_2}, \cdots, \omega  ^{l_s} )$.

 Now we compute the fixed point set $F_{(\phi
,\sigma)}$ of $(\phi,\sigma)$.
If $\left(\begin{array}{c} j\\
g _{\vec{k}_j}\\
\chi _{\vec{l}_j}\end{array}\right)\in F_{(\phi ,\sigma)}$, then
$$\left(\begin{array}{c} \sigma(j)\\
\phi (g _{\vec{k}_j})\\
 \chi _{\vec{l}_j } \phi ^{-1}\end{array}\right)=\left(\begin{array}{c} j\\
g _{\vec{k}_j}\\
\chi _{\vec{l}_j}\end{array}\right)= \left(\begin{array}{c} \sigma (j)\\
g _{\vec{k}_{\sigma (j)}}\\
\chi _{\vec{l}_{\sigma (j)}}\end{array}\right).$$  Thus
\begin {eqnarray} \label {e4.1}A\vec{k}_j=\vec{k}_{\sigma(j)}\ \ \hbox  { and } \ \ \ B\vec{l}_j=\vec{l}_{\sigma(j)},\qquad j=1,\cdots,n. \end
{eqnarray}

 Assume that  $\tau=(j_1,\cdots,j_r)$ is a cycle in independent
 cycle decomposition of   $\sigma$ with length $r$. We have
\begin {eqnarray} \label {e4.2}
\begin{array} {llll}
\vec k_{j_2}&= A  \vec k_{j_1},& \vec l_{j_2}&=B \vec l_{j_1},\\
&\cdots & \cdots & \\
\vec k_{j_r}&=A^{r-1}\vec k_{j_1}, &\vec l_{j_r}&=B^{r-1}\vec l_{j_1},\\
\vec k_{j_1}&=A^r\vec k_{j_1},& \vec l_{j_1}&=B^r\vec l_{j_1}.
\end{array}
\end  {eqnarray}  Thus we only need choice $\vec k_{j_1}$ and $\vec l_{j_1}$
such that $\vec k_{j_1}$ and $\vec l_{j_1}$ become the solutions of
\begin {eqnarray} \label {e4.3}(A^r-I)X=0 \end {eqnarray} and
  \begin {eqnarray} \label {e4.4}(B^r-I)X=0, \ \end {eqnarray} respectively. According to (\ref {e4.2}), we compute
  $\vec
k_{j_2},\cdots,\vec k_{j_r}$ and $\vec l_{j_2},\cdots,\vec
l_{j_r}$. That is, we get $ \vec k_j $' s  and  $ \vec l_j{} $' s
such that (\ref {e4.1}) holds.  Conversely it is true, too.

Let  $\nu _{\phi, r}$ and   $\mu _{\phi, r}$ denote the numbers of
solutions of  (\ref {e4.3}) and  (\ref {e4.4}) in  $ \left (
\begin {array} {l} \mathbb{Z}_{m_1}
\\
\vdots\\
\mathbb{Z}_{m_s} \end {array}  \right)$ and  $ \left (
\begin {array} {l} \mathbb{Z}_{m}
\\
\vdots\\
\mathbb{Z}_{m} \end {array}  \right)$, respectively. If  the type of
$\sigma$ is $1^{\lambda_1}2^{\lambda_2}\cdots n^{\lambda_n}$, then
\begin {eqnarray} \label {e3.5}
|F_{(A,\sigma)}|=\prod_{r=1}^{n} (\nu _{\phi, r}\mu _{\phi,
r})^{\lambda_r}.\end {eqnarray}
\begin{Theorem}\label{3.2}

If  $G$ is a commutative group with order $m$, then
$$
{\cal N}(G, n)=\frac{1}{ \mid \! Aut (G)\! \mid
}\sum\limits_{\phi\in Aut (G)}\
\sum\limits_{\lambda_1+2\lambda_{2}+\cdots +n\lambda_n=n}\
\prod\limits_{i=1}^{n}
  \frac{(\nu _{\phi, i} \mu _{\phi, i}   )^{\lambda_i}}{\lambda_i!i^{\lambda_i}}.
$$
\end{Theorem}

{ \bf Proof.}  By  Lemma \ref {3.1} and  (\ref {e3.5}), we have
\begin{eqnarray*}
{\cal N}(G, n)&=&\frac{\sum\limits_{(\phi,\sigma)\in M}|F_{(\phi,\sigma)}|}{|M |}\nonumber\\
&=&\frac{\sum\limits_{\phi\in Aut (G)}\
\sum\limits_{\lambda_1+\cdots
+n\lambda_n=n}\frac{n!}{\lambda_1!\cdots\lambda_n!1^{\lambda_1}\cdots
n^{\lambda_n}}\ \prod\limits_{r=1}^{n}(\nu _{\phi, r}\mu _{\phi, r})^{\lambda_r}}
{ \mid \! Aut (G) \! \mid  n!}\nonumber\\
&=&\frac{1}{ \mid \! Aut (G) \! \mid }\sum\limits_{\phi\in Aut (G)}
\ \sum\limits_{\lambda_1+2\lambda_{2}+\cdots +n\lambda_n=n}\
\prod\limits_{i=1}^{n}\frac{(\nu _{\phi, i}\mu _{\phi,
i})^{\lambda_i}}{\lambda_i!i^{\lambda_i}}.\Box
\end{eqnarray*}

\renewcommand{\baselinestretch}{1}

%\end{CJK*}

\end{document}